\documentclass[]{article}
\usepackage{tensor}
\usepackage{subcaption}
\usepackage{graphicx}

\usepackage{amsfonts}
\usepackage{amsmath}
\usepackage{amssymb}
\usepackage{mathtools}
\usepackage{stmaryrd}
\usepackage{xstring}

\def\boldn{\textbf{n}}
\def\boldr{\textbf{r}}
\newcommand{\bOmega}{\mathbf{\Omega}}

\usepackage[numbers]{natbib}
\begin{document}
\newcommand{\FixRef}[3][sec:]
{\IfBeginWith{#2}{#3}
        {\StrBehind{#2}{#3}[\RefResult]}
        {\def\RefResult{#2}}\IfBeginWith{#1}{#3}
        {\StrBehind{#1}{#3}[\RefResultb]}
        {\def\RefResultb{#1}}}

\newcommand{\secref}[1]
{\FixRef{#1}{sec:}Section~\ref{sec:\RefResult}}
\newcommand{\secreff}[1]
{\FixRef{#1}{sec:}in Section~\ref{sec:\RefResult}}
\newcommand{\Secreff}[1]
{\FixRef{#1}{sec:}In Section~\ref{sec:\RefResult}}
\newcommand{\secrefm}[2]
{\FixRef[#2]{#1}{sec:}Sections~\ref{sec:\RefResult}-\ref{sec:\RefResultb}}
\newcommand{\secreffm}[2]
{\FixRef[#2]{#1}{sec:}in Sections~\ref{sec:\RefResult}-\ref{sec:\RefResultb}}
\newcommand{\Secreffm}[2]
{\FixRef[#2]{#1}{sec:}In Sections~\ref{sec:\RefResult}-\ref{sec:\RefResultb}}
\newcommand{\figref}[1]
{\FixRef{#1}{fig:}Figure~\ref{fig:\RefResult}}
\newcommand{\figrefm}[2]
{\FixRef[#2]{#1}{fig:}Figures~\ref{fig:\RefResult}-\ref{fig:\RefResultb}}
\newcommand{\figreff}[1]
{\FixRef{#1}{fig:}in Figure~\ref{fig:\RefResult}}
\newcommand{\figreffm}[2
]{\FixRef[#2]{#1}{fig:}in Figures~\ref{fig:\RefResult}-\ref{fig:\RefResultb}}
\newcommand{\Figreff}[1]
{\FixRef{#1}{fig:}In Figure~\ref{fig:\RefResult}}
\newcommand{\Figreffm}[2]
{\FixRef[#2]{#1}{fig:}In Figures~\ref{fig:\RefResult}-\ref{fig:\RefResultb}}
\newcommand{\tabref}[1]
{\FixRef{#1}{tab:}Table~\ref{tab:\RefResult}}
\newcommand{\tabreff}[1]
{\FixRef{#1}{tab:}in Table~\ref{tab:\RefResult}}
\newcommand{\Tabreff}[1]
{\FixRef{#1}{tab:}In Table~\ref{tab:\RefResult}}
\newcommand{\tabrefm}[2]
{\FixRef[#2]{#1}{tab:}Tables~\ref{tab:\RefResult}-\ref{tab:\RefResultb}}
\newcommand{\tabreffm}[2]
{\FixRef[#2]{#1}{tab:}in Tables~\ref{tab:\RefResult}-\ref{tab:\RefResultb}}
\newcommand{\Tabreffm}[2]
{\FixRef[#2]{#1}{tab:}In Tables~\ref{tab:\RefResult}-\ref{tab:\RefResultb}}
\newcommand{\egyref}[1]
{\FixRef{#1}{eq:}Equation~\ref{eq:\RefResult}}
\newcommand{\egyreff}[1]
{\FixRef{#1}{eq:}in Equation~\ref{eq:\RefResult}}
\newcommand{\Egyreff}[1]
{\FixRef{#1}{eq:}In Equation~\ref{eq:\RefResult}}
\newcommand{\egyrefm}[2]
{Equations~\ref{eq:#1}-\ref{eq:#2}}
\newcommand{\egyreffm}[2]
{\FixRef[#2]{#1}{eq:}in Equations~\ref{eq:\RefResult}-\ref{eq:\RefResultb}}
\newcommand{\Egyreffm}[2]
{\FixRef[#2]{#1}{eq:}In Equations~\ref{eq:\RefResult}-\ref{eq:\RefResultb}}
\newcommand{\charef}[1]
{\FixRef{#1}{cha:}Chapter~\ref{cha:\RefResult}}
\newcommand{\chareff}[1]
{\FixRef{#1}{cha:}in Chapter~\ref{cha:\RefResult}}
\newcommand{\Chareff}[1]
{\FixRef{#1}{cha:}In Chapter~\ref{cha:\RefResult}}
\title{An Angular Multigrid Preconditioner for the Radiation Transport Equation with Forward-Peaked Scatter}

\author{Danny Lathouwers and Zolt\'an Perk\'o \\
        d.lathouwers@tudelft.nl , z.perko@tudelft.nl \\
        Department of Radiation Science and Technology, \\
        Delft University of Technology, \\
        Mekelweg 15, 2629JB Delft, The Netherlands}

\maketitle

\begin{abstract}
In a previous paper (Lathouwers and Perk\'o, 2019) we have developed an efficient angular multigrid preconditioner
for the Boltzmann transport equation with forward-peaked scatter modeled by the Fokker-Planck approximation. The
discretization was based on a completely discontinuous Galerkin finite element scheme both for space and angle.
The scheme was found to be highly effective on isotropically and anisotropically refined angular meshes.
The purpose of this paper is to extend the method to non-Fokker-Planck models describing the forward scatter by 
general Legendre expansions. As smoother the standard source iteration is used whereas solution on the coarsest angular mesh is
effected by a special sweep procedure that is able to solve this problem with highly anisotropic scatter using only a small number of
iterations. An efficient scheme is obtained by lowering the scatter order in the multigrid preconditioner. A set of test
problems is presented to illustrate the effectivity of the method, i.e. in less iterations than the single-mesh case
and more importantly with reduced computational effort.
\end{abstract}
%
%
%
\section{Introduction}
\label{sec:Introduction}
The solution of the linear Boltzmann equation (LBE) is computationally expensive due to its high dimensionality. Studies of 
numerical techniques focus mostly on the steady mono-energetic form being the building block for more complex cases.
Discretization is generally performed by discrete ordinates which is straightforward to implement. 
Spatial discretization is most often performed by the discontinous Galerkin in space. Whatever the 
discretization used, some form of iteration is needed to resolve the system of equations.
Source iteration, the most widespread technique, works well for optically thin media but is less suitable for thick
scattering materials. As the slowest mode of convergence in thick scattering media has little angular dependence this can be remedied by
the diffusion synthetic acceleration (DSA) that accelerates reduction of these errors by solution of a diffusion 
problem. When used as preconditioner in a Krylov method, this practically gives rise to an unconditionally effective and stable method.

For anisotropic scattering and in thin media, DSA is still ineffective. Multigrid methods where the problem is formulated
on multiple levels lead to a more successfull approach. The basic principle of the multigrid approach is to smooth small wavelength
errors on successive grids combined with an exact solver on the coarsest grid. DSA can be viewed as a two-grid technique in physical
approach, switching between complete transport on the fine level and diffusion on the coarse level.

Our previous paper \cite{Lathouwers2019} has an extensive discussion of literature on multigrid approaches which we only briefly summarize here.
The multigrid technique for radiation transport problems was pioneered by Morel and Manteuffel \cite{Morel1991} where they 
constructed an angular multigrid method for the one dimensional
$S_N$ equations. The method was extended to two dimensions by Pautz \cite{Pautz1999} by the introduction of high-frequency filtering to increase
stability of the method. Various references have investigated the efficiency of angular multigrid schemes for both the $P_N$ and the $S_N$ equations
\cite{Lee2010a,Lee2010b,Oliveira2000,Turcksin2012,Drumm2017}. These papers have shown the great benefit of multigrid: speedups of up to 10 are reported 
compared to single grid methods.

In previous work we have introduced a novel finite element angular discretization of the transport equation \cite{Kophazi2015}. 
The method offers the capability to anisotropically refine the sphere to focus on important directions.
In a later paper we added a discretization
scheme for the Fokker-Planck small angle scatter term \cite{Hennink2017} and an angular multigrid scheme for the efficient solution
\cite{Lathouwers2019}. The scheme utilizes the multigrid method as preconditioner for a Krylov method and was found to be highly effective 
compared to the single mesh preconditioner. The purpose of this paper is to widen the scope of these previously introduced
methods to the more general case of highly anisotropic scatter as modeled by Legendre scatter expansion.

This paper is outlined as follows. In Section 2, we summarize the discontinuous Galerkin discretization in space-angle that our
method is based on. The standard source iteration procedure is described in Section 3. The angular multigrid procedure proposed
in the present work is discussed in Section 4 comprising of the mesh hierarchy, and the inter-grid transfers. 
Section 5 describes a newly
formulated sweep methodology for efficient coarse mesh solution compatible with high order scatter. The method is tested on a set of model problems. Final conclusions are drawn in Section 6.
\section{Space-angle discretization of the transport equation}
\label{sec:Space_angle_discretization}
The space-angle discretization of the Boltzmann equation has been presented in detail in \cite{Kophazi2015} and \cite{Hennink2017}, and will 
therefore be described briefly by only covering the essentials. Details can be found in the original references.
Particle transport with highly anisotropic scatter is described by the linear Boltzmann transport equation.
We neglect energy dependence in this paper as the focus is on accelerating the single group iterative method.
The linear mono-energetic Boltzmann equation reads
\begin{equation}
  \bOmega \cdot \nabla \phi(\boldr,\bOmega) + \Sigma_t(\boldr) \phi(\boldr,\bOmega) = 
    Q(\boldr,\bOmega)
  \label{eq:LBE}
\end{equation}
where $\boldr$ is the spatial coordinate, $\bOmega$ is the unit direction vector, $\phi$ is the angular flux density, $Q$ is the 
volumetric source density including scatter, $\Sigma_t$ is the total macroscopic cross section. 
The angular flux for incoming directions is specified on the domain boundary, $\Gamma_I$,
\begin{equation}
  \phi (\boldr,\bOmega) = \phi_{in} (\boldr,\bOmega), \; \boldr \in \Gamma_I, \; \bOmega \cdot \hat{\boldn} < 0
\end{equation}
\subsection{Phase space elements}
The spatial domain $V$ is made up of elements $V_k$, where $k$ is the index of the spatial element.
A discontinuous solution space $S_{h,p}$, is defined containing polynomials of order $p$ at most. This is a
standard approach and we focus our attention on the angular discretization.

The construction of angular elements is based on hierarchical sectioning of the unit sphere into patches, $D_p$,
where $p$ is the patch index. The coordinate planes divide the sphere 
into cardinal octants, which are spherical triangles. We also assign a level, $l_p$, to a patch.
The spherical triangles at the coarsest level are assigned a level of $l_p = 0$.
Spherical triangles with $l_p = 1$ are obtained by halving the edges of the $l_p = 0$ patches and subsequently 
connecting the emerged points with
great circles. Every patch is hereby split into four daughters. This procedure can be repeated to arbitrary depth and 
locally on the sphere (aniostropic refinement).

The angular subdivision of the sphere is described by a set $P$ of patch indices such that
$\cup_{p \in P} D_p \equiv D$ where $D$ denotes the unit sphere surface and $D_p \cap D_q \equiv \emptyset$ for $\forall p,q \in P, p \neq q$. The phase 
space mesh is then obtained by assigning an angular subdivision $P_k$ to each spatial element $V_k$ providing a high level of flexibility.
\subsection{Angular basis functions}
Two sets of angular basis functions, $\psi_{[p] d} (\bOmega)$, are used throughout this paper. Both are local to the patch $D_p$ by setting
$\psi_{[p] d} (\bOmega)=0$ if $\bOmega \notin D_p$ and are discontinuous at the patch boundary.
Here $\psi_{[p] d} (\bOmega)$ denotes the $d$-th basis function on the patch with index $p$. 
The locality-property of the functions ensures that the streaming-removal terms will not couple to non-overlapping patches.
As explained later, this eases the use of sweep-based algorithm.
The two basis functions sets are:
\begin{enumerate}
  \item{\textbf{Const} A unit value function on the patch.}
  \item{\textbf{Lin} A nodal set of three functions satisfying $\psi_{[p] d} (\bOmega_{d'}) = \delta_{dd'}$
        where $\bOmega_{d'}$ are the patch-vertices. These functions result from projecting the standard 
        Lagrange functions on a specific flat triangle on the octahedron onto the sphere. The flat triangle is formed by 
        projecting the vertices of the patch to the octahedron.}
\end{enumerate}
\subsection{Discretization}
The flux in each energy group is written as a product of spatial and angular basis functions as
\begin{equation}
  \phi (\boldr,\bOmega) = \sum_{k,i} \sum_{p \in P_k,d} \phi\indices{_{[k,p]}^{i,d}} \Phi_{[k]i} (\boldr) \Psi_{[p]d} (\bOmega)
  \label{eq:Expansion}
\end{equation}
where $\Phi_{[k]i} (\boldr)$ is the $i$-th spatial basis function of element $k$ and $\Psi_{[p]d} (\bOmega)$ is the $d$-th
angular basis function on patch (angular element) $p$. 

Substituting the expansion given by \egyref{eq:Expansion} into the transport equation (\egyref{LBE}), multiplying with a test function in space and angle 
and subsequently integrating over complete phase space leads to:
\begin{equation}
  \label{eq:complete}
  \sum_f \Upsilon_{[f,j,q]lm} - 
  \sum_{\xi=1}^3 \sum_i \sum_d V_{[j]li\xi} \phi\indices{_{[j,q]}^{id}} A\indices{_{[q,p]md}^\xi} +
  \Sigma_a \sum_i \sum_d N_{[j]li} \phi\indices{_{[j,q]}^{id}} M_{[q]md} = 
  Q_{[j,q]lm},
\end{equation}
where the streaming term reads
\begin{equation}
  \Upsilon_{[f,j,q]lm} = \int_{\partial V_j^f} \Phi_{[j]l} (\boldr) 
  \sum_{\begin{array}{c} k \in \{j,j_f'\} \\ i \end{array}} 
  \sum_{\begin{array}{c} p \in P_k \\ d \end{array}} 
  \phi\indices{_{[k,p]}^{id}} \Phi_{[k]i} (\boldr) 
  A_{[f,q,p]md} d \boldr
  \label{eq:Surface_Integral}
\end{equation}
and we made use of the following shorthand notations:
\begin{eqnarray}
  N_{[j]li} = \int_{V_j} \Phi_{[j]l} (\boldr) \Phi_{[j]i} (\boldr) d \boldr \\
  V_{[j]li\xi} = \int_{V_j} \nabla_\xi \Phi_{[j]l} (\boldr) \Phi_{[j]i} (\boldr) d \boldr \\
  A\indices{_{[q,p]md}^\xi} = \int_{D_q} \Omega^\xi \Psi_{[q]m} (\bOmega) \Psi_{[p]d} (\bOmega) d \bOmega \\
  A_{[f,q,p]md} = \sum_{\xi=1}^3 {\hat{n}}_{[f] \xi} A\indices{_{[q,p]md}^{\xi}} \\
  A_{[f,q]md} = A_{[f,q,q]md} \\
  M_{[q,p]md} = \int_{D_q} \Psi_{[q]m} (\bOmega) \Psi_{[q]d} (\bOmega) d \bOmega \\
  M_{[q]md} = M_{[q,q]md} \\
  Q_{[j,q]lm} = \int_{V_j} \int_{D_q} \Phi_{[j]l} (\boldr) \Psi_{[q]m} (\bOmega) Q(\boldr,\bOmega) d \bOmega d \boldr.
\end{eqnarray}
Here, element $j$ has faces indexed by $f$ with its neighbor at face $f$ denoted as $j_f'$ and ${\hat{n}}_{[f] \xi}$ is 
the $\xi$ coordinate of the outward normal of face $f$ in element $j$. The mass matrices in space and angle are denoted by $N$ and $M$, respectively.
$V$ contains the volmetric streaming integral and $A$ is the angular Jacobian.

The face integrals from streaming need to be made unique by an upwinding procedure. In $S_N$ schemes this is straightforward.
Here, the angular components on a patch need to be treated simultaneously.
Various authors (e.g. \cite{Pain2006}) have introduced Riemann procedures to separate the surface terms into inward and outward contributions.
In previous papers \cite{Kophazi2015,Hennink2017} we have derived how the numerical flux needs to be evaluated in the different situations, i.e. where
the neighbor is either (i) equally refined, (ii) coarser or (iii) finer. To prevent repetition and concentrate
on the multigrid aspects, we refer the reader to the original references for more in-depth discussion.
\section{Single-grid solution approach}
\label{sec:solver}
The discrete transport equation is written as
\begin{equation}
  L {\boldsymbol \phi} = S {\boldsymbol \phi} + {\bf f}.
\end{equation}
Here, $L$ and $S$ are the discrete transport and scatter operators, respectively and $f$ is
the vector containing the independent source.

We use preconditioned Bicgstab \cite{bicgstab} to solve the linear system with
$L^{-1}$ as preconditioner, corresponding to Krylov-accelerated source iteration
\begin{equation}
  L {\boldsymbol \phi}^{k+1} = S {\boldsymbol \phi}^k + {\boldsymbol f}.
\end{equation}
In discrete ordinates codes the operator $L^{-1}$ is easily applied as the ordinates are independent.
Here, we use a finite element discretization in angle, where such a sweep is no longer exact. The cause is possible
bi-directionality due to some directions on a patch being incoming whereas others are outgoing with respect to
the face.
To precondition the system we have devised the following sweep algorithm:
\begin{itemize}
	\item An $S_2$ ordinate set is used and for each ordinate, corresponding to a particular octant, the sweep order of the 
              elements is determined.
	\item For each direction, we traverse the spatial elements and each 
	angular element in the octant is visited and the corresponding linear system is solved.
	In more mathematical terms we perform
	a block Gauss-Seidel iteration for a given ordering of the spatio-angular elements. Splitting $L$ 
	into implicit and explicit part, $L=L_I + L_E$, the iteration reads
	\begin{equation}
	L_I {\bf \phi}^{k+1} = (S - L_E) {\bf \phi}^k + {\bf f}
	\end{equation}
	and the preconditioner then is $L_I$.
\end{itemize}
We have previously demonstrated that this sweep algorithm is an effective preconditioner (see \cite{Kophazi2015} for more details). 
In problems where scatter is highly anisotropic the procedure unfortunately becomes less effective.
\section{An angular multigrid preconditioner}
\label{sec:Angular_Multigrid_Preconditioner}
The multigrid method smooths the error on a given mesh
and then transfers the residual to a coarser mesh where the lower frequency errors can be
more effectively attenuated. As we deal with a linear problem, we use the linear multigrid algorithm as 
shown in Figure~\ref{fig:lmg}. 
\begin{figure}
  Algorithm $LMG(\boldsymbol{\phi}_l,{\bf{f}}_l,l)$ \\
  if $(l=0)$ then \\
    \mbox{} \mbox{} $S(\boldsymbol{\phi}_l,{\bf{f}}_l,\nu_{coarse})$ \\
  else \\
    \mbox{} \mbox{} $S(\boldsymbol{\phi},f,\nu_{pre})$ \\
    \mbox{} \mbox{} ${\bf{r}}_l = {\bf{f}}_l - A_l \boldsymbol{\phi}_l$ \\
    \mbox{} \mbox{} ${\bf{f}}_{l-1} = R {\bf{r}}_l$ \\
    \mbox{} \mbox{} $\boldsymbol{\phi}_{l-1} = 0$ \\
    \mbox{} \mbox{} $LMG(\boldsymbol{\phi}_{l-1},{\bf{f}}_{l-1},l-1)$ \\
    \mbox{} \mbox{} $\boldsymbol{\phi}_l = \boldsymbol{\phi}_l + P \boldsymbol{\phi}_{l-1}$ \\
    \mbox{} \mbox{} $S(\boldsymbol{\phi}_l,{\bf{f}}_l,\nu_{post})$ \\
  endif \\
  end Algorithm LMG \\
  \caption{Recursive linear multigrid algorithm based on the V-cycle (\cite{Wesseling1992}). The number of pre- and post-smoothing 
           steps is $\nu_{pre}$ and $\nu_{post}$,
           respectively. The restriction and prolongation operators are denoted by R and P.}
  \label{fig:lmg}
\end{figure}
Since the multigrid method has been extensively documented (ee e.g. \cite{Wesseling1992}), we only discuss the particular multigrid 
components.
\subsection{Multigrid components}
\begin{figure}[!ht]
	\centering
	\begin{subfigure}[b]{0.32\textwidth}
	\includegraphics[width=0.9\columnwidth]{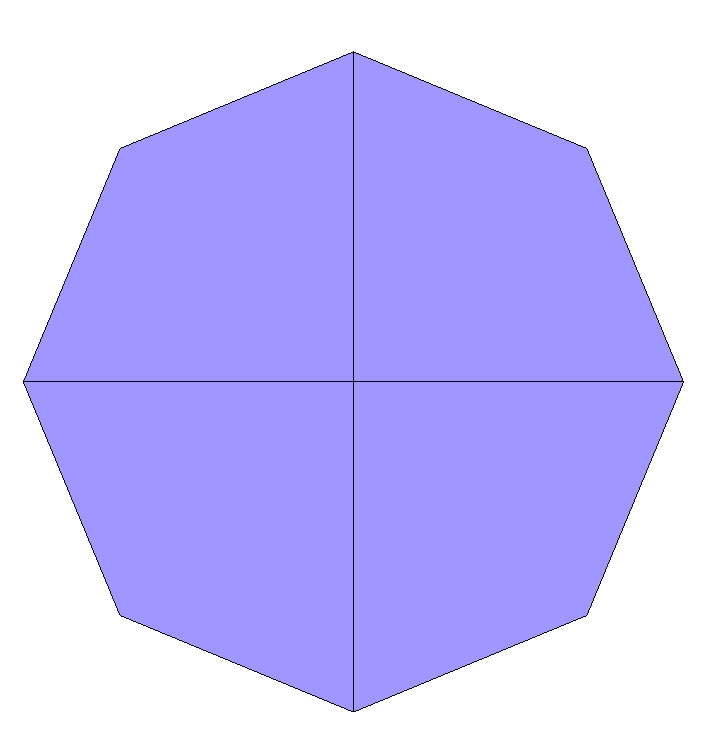}
	\caption{$T_0$}
	\end{subfigure}
	\begin{subfigure}[b]{0.32\textwidth}
	\includegraphics[width=0.9\columnwidth]{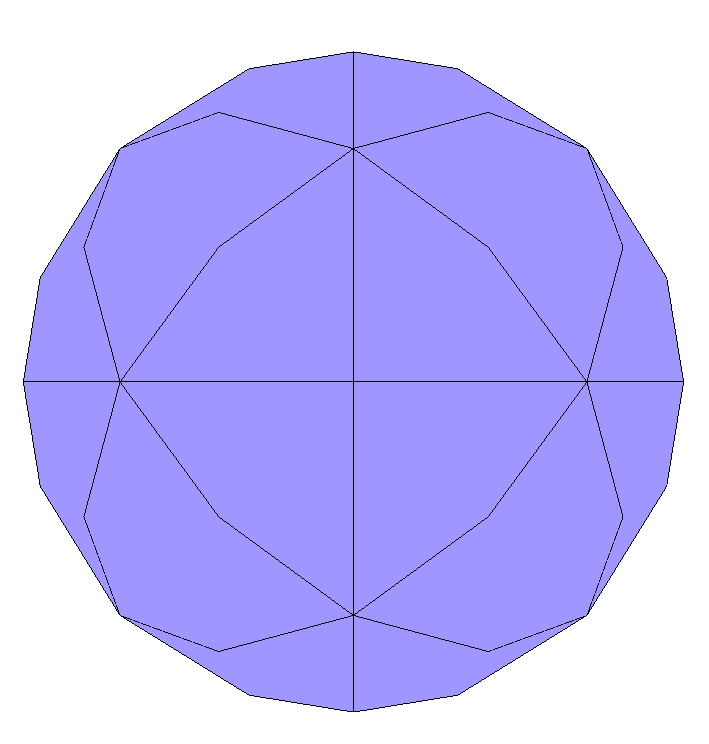}
	\caption{$T_1$}
	\end{subfigure}
	\begin{subfigure}[b]{0.32\textwidth}
	\includegraphics[width=0.9\columnwidth]{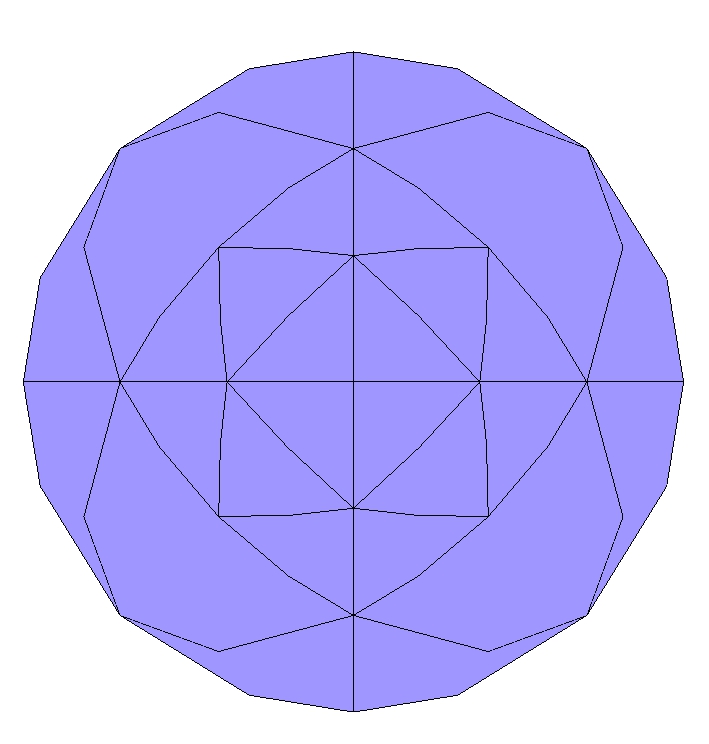}
	\caption{$T_2$}
	\end{subfigure} \\
	\begin{subfigure}[b]{0.32\textwidth}
	\includegraphics[width=0.9\columnwidth]{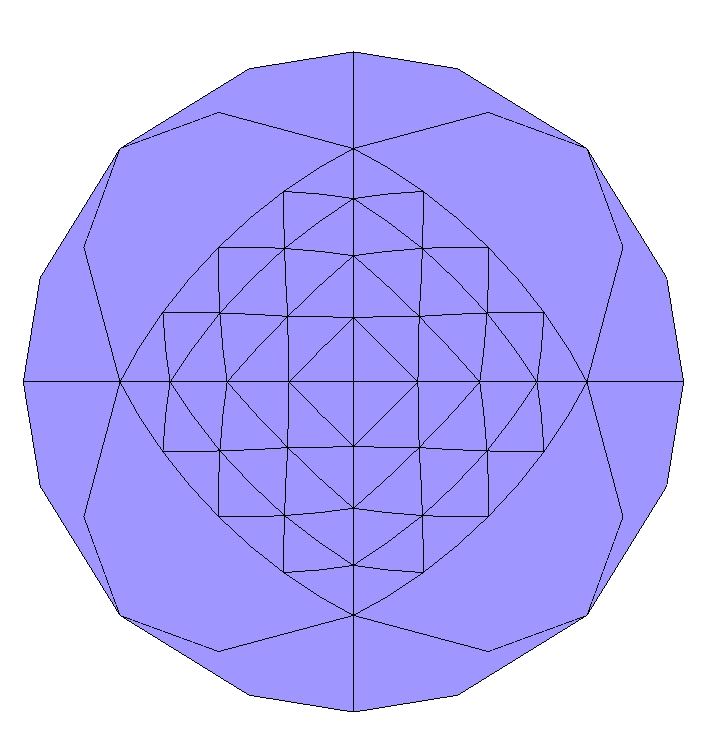}
	\caption{$T_3$}
	\end{subfigure}
	\begin{subfigure}[b]{0.32\textwidth}
	\includegraphics[width=0.9\columnwidth]{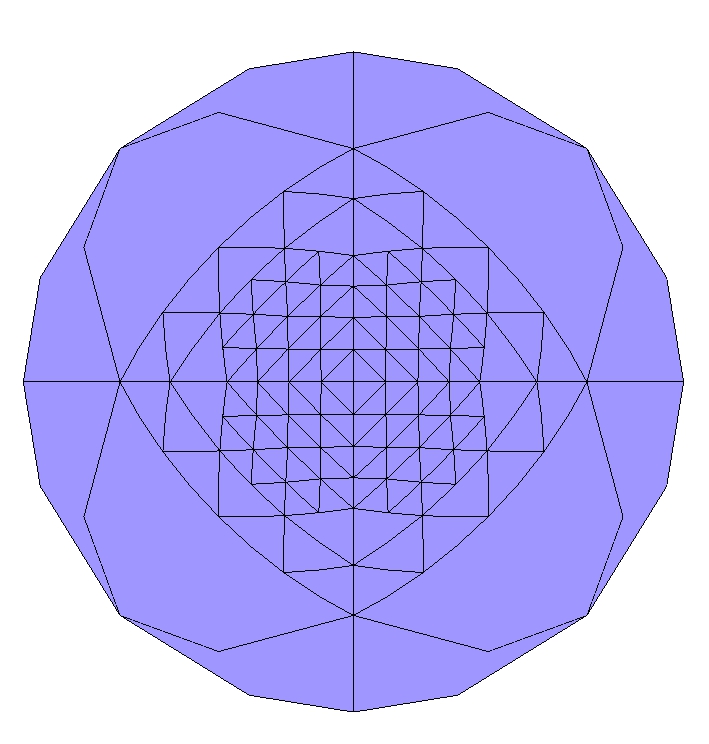}
	\caption{$T_4$}
	\end{subfigure}
	\begin{subfigure}[b]{0.32\textwidth}
	\includegraphics[width=0.9\columnwidth]{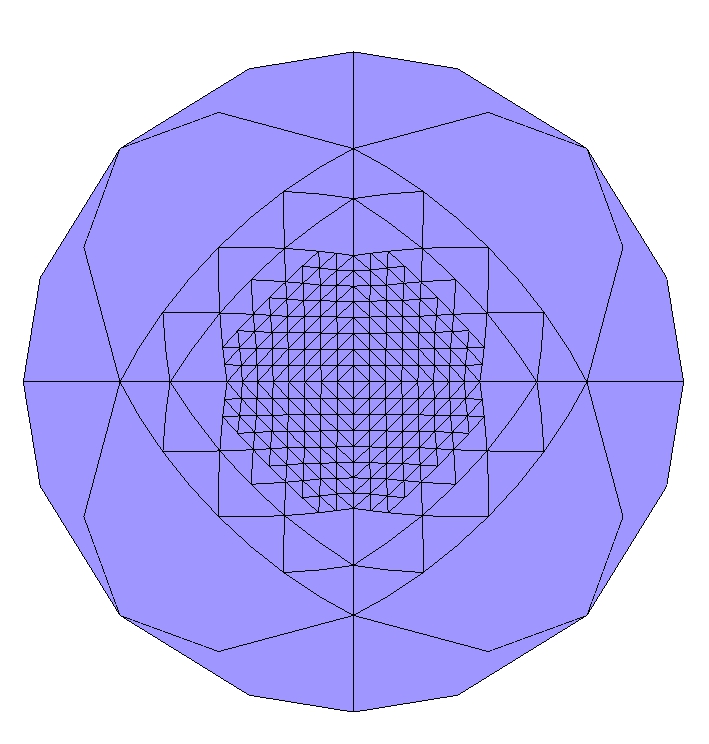}
	\caption{$T_5$}
	\end{subfigure}
	\caption{A set of nested angular meshes. Mesh $T_5$ is the finest mesh. The mesh elements 
                 are shown as bisected triangles rather than as spherical triangles. This is due to the plotting.}
	\label{fig:circ_meshes}
\end{figure}

A nested series of spherical meshes is obtained by refining a uniformly discretized sphere consisting of 8 spherical 
triangles (the octants).
The coarsest mesh, $T_0$ has level $0$. Refined meshes are constructed by refining specific angular elements.
Angular meshes are constrained to possess only
up to two-irregularity, i.e. neighboring angular elements can differ by two levels at most.
A series of triangulations $\{T_l\}$ is obtained with the maximum angular element level on $T_l$ is $l_p=l$.
The finest angular mesh is $T_L$. The maximum occuring element refinement in the problem is $l_p=L$.
An example of a series of spherical meshes is shown \figreff{fig:circ_meshes}.

The prolongation operator follows naturally from the use of a nested discontinuous finite element space, i.e. the complete angular solution 
can be prolongated without approximation from the coarse to the fine grid by Galerkin projection. The restriction operator is the 
transpose of prolongation operator, i.e. $R = P^T$.

As smoother the standard source iteration procedure is used. Hence we perform iterations of the form
\begin{equation}
  L_l \phi_l^{k+1} = S_l \phi_l^k + f_l
\end{equation}
where k is the iteration index.

The coarse mesh problems are formulated by direct discretization (Discretization Coarse Grid Approximation, DCA). 
This procedure is compatible with our matrix-free implementation of the transport solver.
Solution on the coarsest level could be performed by source iteration. This is however not effective for highly anisotropic scatter. 
In Section~\ref{subsec:alt_sweep} we describe a more effective approach.
\section{Results}
\label{sec:Results}
\begin{figure}
  \includegraphics[width=0.9\columnwidth]{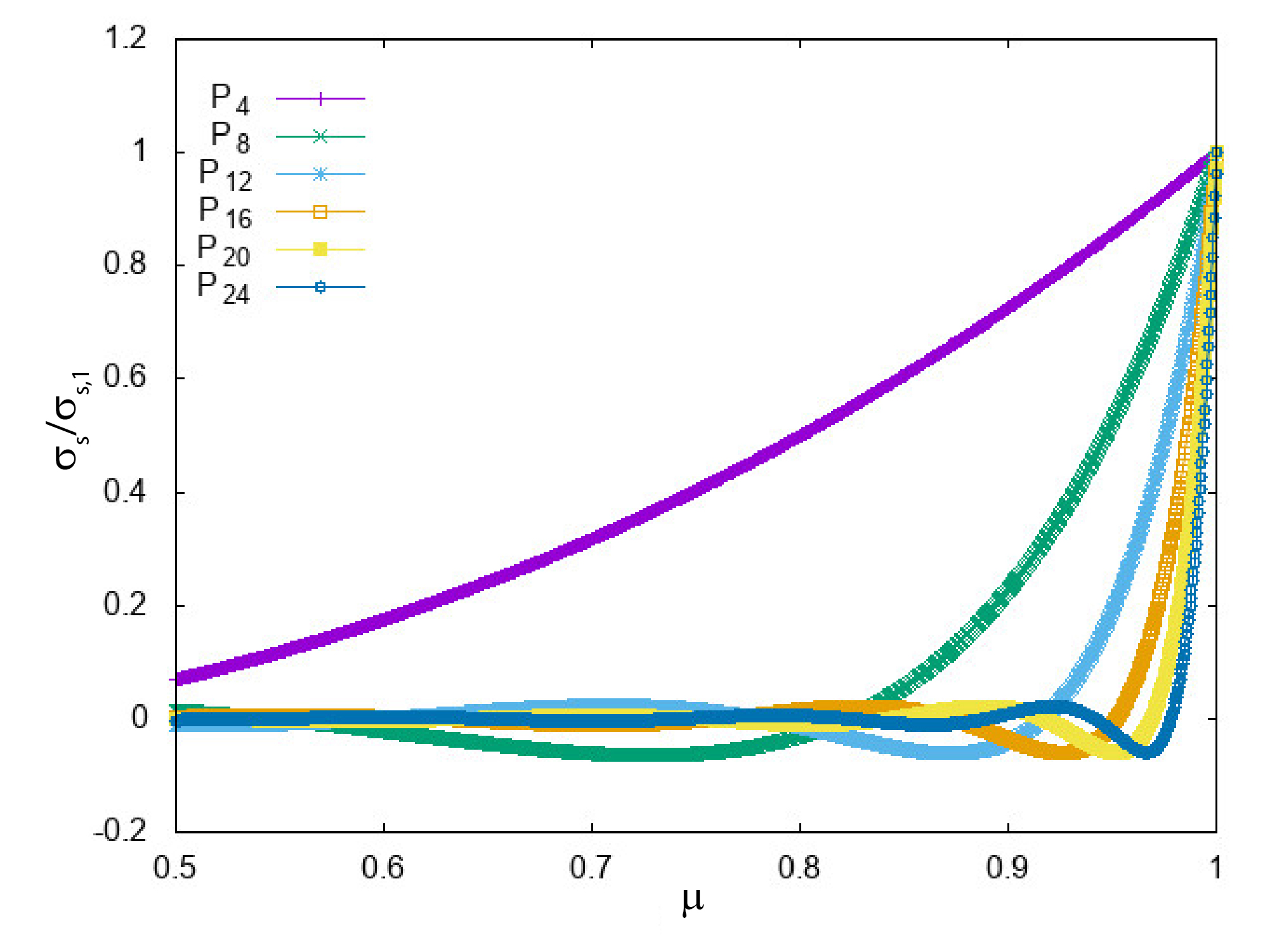}
  \caption{Scatter cross section normalized to maximum value as function of $\mu$ for different expansion orders, $N$.}
  \label{fig:scatter_for_diff_orders}
\end{figure}

Although our method is general in terms of scattering functions used we focus on the Fokker-Planck equivalent scatter kernel.
For such scatter kernels we have that the Legendre moments are given by
\begin{equation}
  \sigma_{s,n} = \frac{\alpha}{2} \left( N(N+1) - n(n+1) \right)
  \label{eqn:scatter_moments}
\end{equation}
where $\alpha$ is the momentum transfer coefficient. In the present paper we will vary the scatter order, $N$, to increase the
forward-peaked nature of the scatter kernel. The normalized scatter cross section for varying $N$ is shown in 
Figure~\ref{fig:scatter_for_diff_orders}.
\subsection{Preliminary multigrid performance}
\label{subsec:prel_mg_perf}
In this section we present a preliminary study of multigrid performance and a discussion of the aspects that are crucial for 
obtaining an efficient method. 

We apply the basic multigrid algorithm to a radiation transport problem in a cube of size $5 \times 5 \times 5\ cm^3$. The
geometry is meshed with $30 \times 30 \times 30$ hexahedral elements. The meshes were produced by Gmsh \cite{gmsh2009}. Boundary 
conditions are vacuum on all sides and the problem is driven by a uniform isotropic unit strength source 
in the cube. The transport cross section is equal to 1 and the scatter moments are calculated from equation~\ref{eqn:scatter_moments}.

We use the multigrid V-cycle with the transport sweep as smoother, inter grid transfers based on Galerkin projection,
discussed previously and a thorough coarse grid solution.  
The stop criterion is chosen as $||r|| / ||b|| < 10^{-8}$ where $r$ is the residual vector and $b$ the right hand side.
The stop criterion for the coarse grid solver is set to $10^{-5}$.
Uniformly refined angular meshes are used. Angular discretization is based on constant and linear basis functions.

For different angular refinement levels on the finest mesh (l = 1, 2, or 3), 
we investigate the effects of the number of smoothing cycles and restricting
the scatter order on the coarse level to $N_r < N$ on the multigrid performance. The results for the number of iterations obtained are given in Table~\ref{tab:opt_mg_perf}.
\begin{table}
\centering
\caption{Iteration counts for isotropically refined angular meshes (l) up to level 3 and scatter orders N of 4,8 and 16 for the 3D box problem with uniform isotropic 
source for single grid (SG), different multigrid cycles and limited scatter orders $N_r$ on the coarsest level.}
\begin{tabular}{|c|c|c|c|c|c|c|c|c|c|c|c|}
\hline
basis & l & $N_r$ &  SG  & V(1,1) & V(2,1)       &  SG  &  V(1,1) & V(2,1)      & SG  &  V(1,1) & V(2,1) \\ \hline \hline
      &   &       & \multicolumn{3}{|c|} {$N=4$} & \multicolumn{3}{|c|} {$N=8$} & \multicolumn{3}{|c|} {$N=16$} \\ \hline \hline
const & 1 &   0   &  39  &  14    &   14         &  93  &    22   &   22        & 179 &    42   & 42  \\
const & 1 &   1   &      &   7    &    7         &      &    12   &   11        &     &    18   & 16  \\
const & 1 &   2   &      &   7    &    7         &      &    11   &   11        &     &    19   & 19  \\ \hline
const & 2 &   0   &  43  &  15    &   14         & 103  &    33   &   33        & 187 &    49   & 42  \\
const & 2 &   1   &      &   9    &    8         &      &    17   &   15        &     &    24   & 23  \\
const & 2 &   2   &      &   9    &    8         &      &    17   &   15        &     &    26   & 23  \\ \hline
const & 3 &   0   &  44  &  17    &   15         & 121  &    34   &   34        & 253 &   111   & 210  \\
const & 3 &   1   &      &   9    &    8         &      &    19   &   18        &     &    50   & 34  \\
const & 3 &   2   &      &   9    &    8         &      &    20   &   17        &     &    40   & 36  \\
const & 3 &   3   &      &        &              &      &    17   &   14        &     &    32   & 28  \\ \hline
lin   & 1 &   0   &  44  &  15    &   14         & 127  &    45   &   33        & 287 &    245  & 231 \\
lin   & 1 &   1   &      &   6    &    5         &      &    15   &   14        &     &     51  & 47  \\
lin   & 1 &   2   &      &   5    &    4         &      &    12   &    9        &     &     26  & 22  \\
lin   & 1 &   3   &      &        &              &      &     9   &    9        &     &     22  & 21  \\ \hline
lin   & 2 &   0   &  45  &  14    &   14         & 123  &    42   &   32        & 294 &         &     \\
lin   & 2 &   1   &      &   5    &    4         &      &    14   &   11        &     &     43  & 41  \\
lin   & 2 &   2   &      &   4    &    4         &      &     9   &    8        &     &     23  & 21  \\
lin   & 2 &   3   &      &        &              &      &     9   &    8        &     &     22  & 18  \\ \hline
lin   & 3 &   0   &  44  &  13    &   12         & 124  &    31   &   31        & 284 &         &     \\
lin   & 3 &   1   &      &   4    &    4         &      &    11   &   12        &     &     31  & 36  \\
lin   & 3 &   2   &      &   4    &    3         &      &     8   &    7        &     &     22  & 22  \\
lin   & 3 &   3   &      &        &              &      &     4   &    3        &     &     22  & 20  \\ \hline
\end{tabular}
\label{tab:opt_mg_perf}
\end{table}
The main conclusions that can be drawn from this table are as follows: 
(i) Source iteration is an 
excellent smoother with V(1,1) not performing significantly worse than V(2,1). Considering the amount of work
in smoothing, V(1,1) is favorable in all cases and V(2,1) is not considered throughout the remainder of the paper. 
(ii) Multigrid performance depends strongly on the choice of scatter order, $N_r$, in the coarse grid problem. The 
iteration counts decrease with increasing $N_r$. The iteration counts are consistently
much lower than when using the single grid solver and more so for higher scatter anisotropy. Timings are not given in 
this table as the coarse grid problem has been deliberately solved accurately to investigate optimal multigrid 
performance and this part dominates the solution time in this approach.
It is concluded that in order to build an effective acceleration technique for anisotropic scatter using multigrid, 
the coarse grid problem needs to be based on a scatter order that is sufficiently high. Using standard source 
iteration as solver on the coarse angular mesh however requires many steps
to reduce the error adequately defying the goal of a cheap solution on the coarse angular mesh. An alternative 
approach is required for the coarse level.
\subsection{Alternative coarse grid solver}
\label{subsec:alt_sweep}
The results using the exact solution at the coarsest grid indicate clearly that the coarse grid operator needs to include sufficiently
high scatter orders for the multigrid algorithm to be efficient. The use of isotropic or linear anisotropic scatter leads 
to poor convergence
ruling out DSA and its variations. Standard source iteration is a (reasonably) efficient smoother but not a good solution method for 
anisotropic scattering media, even on the coarsest angular grid.

Our previous work on the Fokker-Planck equation \cite{Lathouwers2019} has shown that for that case a good solver is obtained by performing a 
block Gauss-Seidel iteration over the angular elements, with 10 such sweeps being sufficient on the coarse level. Here, in the Legendre-scatter case we follow the same route: instead of using source iteration, a Gauss-Seidel technique is used.
Such a Gauss-Seidel technique is easily implemented by first expanding the source vector given by the usual Legendre summation
\begin{equation}
  Q(\boldr,\bOmega) = \sum_{n,o} \sigma_{s,n} \Phi_{n,o} (\boldr) Y_{n,o}(\bOmega)
\end{equation}
with the flux moments defined as
\begin{equation}
  \Phi_{n,o} (\boldr) = \int_{4 \pi} \phi(\boldr,\bOmega) Y_{n,o}(\bOmega)d \bOmega
\end{equation}
Inserting and expanding the flux moments leads to the discrete source, $Q_{[j,q]lm}$
\begin{equation}
  Q_{[j,q]lm} = \sum_{p \in P_j,d,i} N_{[j]li} \sum_{n,o} \sigma_{s,n} X_{[p],d,n,o} X_{[q],m,n,o}\phi\indices{_{[j,p]}^{id}}
\end{equation}
where
\begin{equation}
  X_{[q],m,n,o} = \int_{D_q} \Psi_{[q]m} (\bOmega) Y_{n,o}(\bOmega) d \bOmega
  \label{eq:X}
\end{equation}
These integrals are pre-calculated and stored. On a given spatial element $j$ and angular element, $q$, a linear system can be 
solved considering all flux unknowns outside $j$,$q$ as given. The order in which the spatial and angular elements are visited is equal to the order in the standard sweep. This iteration has shown to be effective on the coarse grid. As in our previous 
work, 10 of such sweeps are sufficient. This coarse grid solver is used throughout the remainder of this work.
\subsection{Reduced multigrid scatter order}
\label{subsec:red_scat_order}
\begin{figure}[!ht]
       \centering
       \begin{subfigure}[b]{0.49\textwidth}
       \includegraphics[width=0.9\columnwidth]{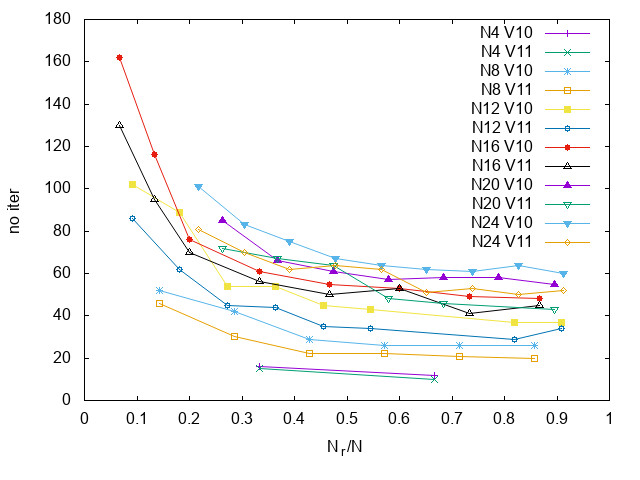}
       \caption{Iteration counts for constant angular basis functions.}
       \end{subfigure}
       \begin{subfigure}[b]{0.49\textwidth}
       \includegraphics[width=0.9\columnwidth]{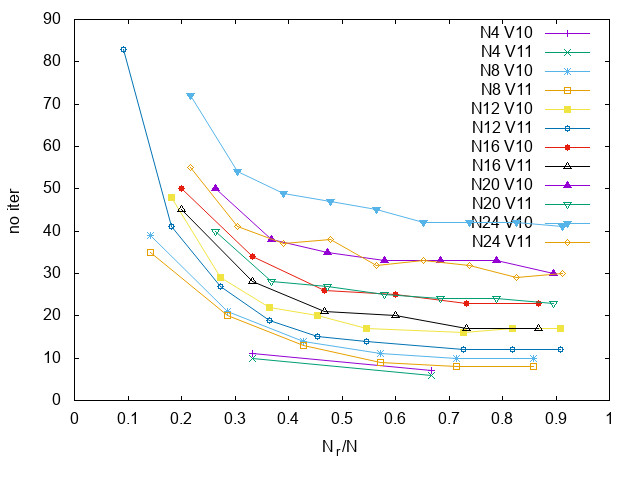}
       \caption{Iteration counts for linear angular basis functions.}
       \end{subfigure} \\
       \begin{subfigure}[b]{0.49\textwidth}
       \includegraphics[width=0.9\columnwidth]{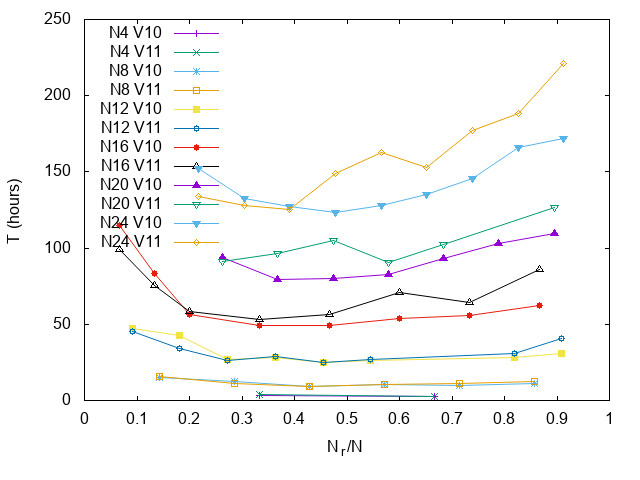}
       \caption{Computational time for constant angular basis functions.}
       \end{subfigure}
       \begin{subfigure}[b]{0.49\textwidth}
       \includegraphics[width=0.9\columnwidth]{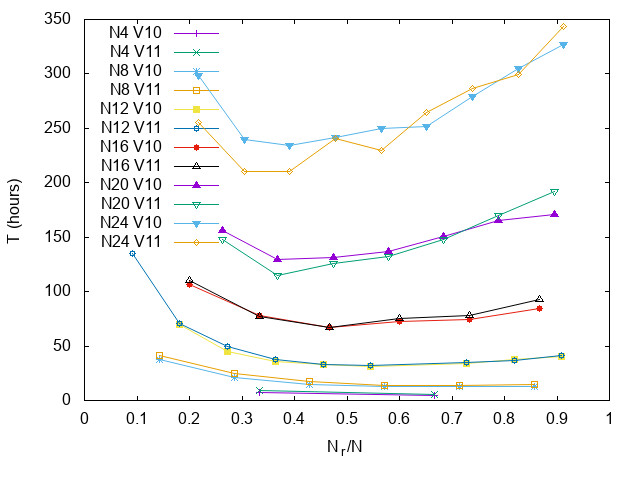}
       \caption{Computational time for linear angular basis functions.}
       \end{subfigure}
       \caption{Effectivity of the multigrid algorithm for varying levels of the scatter order for the multigrid preconditioner using the V(1,0) and V(1,1) cycles for the 3D box geometry with uniform source and isotropic angular refinement. Iteration counts given in a and b. Computational times are given in c and d. Constant angular basis functions used in a and c. Linear basis is used in b and d.}
       \label{fig:Lr_effect}
\end{figure}

The algorithm discussed in Section~\ref{subsec:alt_sweep} - replacing the coarse mesh solver by the alternative sweep solver - is not optimal. Especially
the large number of scatter moments involved for large $N$ is detrimental for performance. In the present section we investigate the 
performance of the multigrid algorithm when operating with reduced scatter order in the multigrid preconditioning stage.
We consider the same problem as before, i.e. the box geometry with isotropic uniform source with an angular
mesh that is 3 times isotropically refined. The multigrid preconditioner 
uses a reduced maximum scatter order $N_r$. We use 10 sweeps on the coarsest grid. The stop criterion is $||r|| / ||b|| < 10^{-8}$.
The results of varying $N_r$ for different angular basis functions, and different cycles (V(1,0) and V(1,1)), are 
shown in Figure~\ref{fig:Lr_effect}.
The number of iterations required as function of $N_r$ decreases with the observation that linear basis function require less
iterations than constant basis functions.
This may be related to the interpolation operators in the constant 
basis function case not being accurate enough: Mesh independent multigrid convergence can only be obtained when the following condition
is met: $m_P + m_R > 2m$ where $m_P$ and $m_R$ are the order of polynomial that can be exactly prolongated, and restricted and $2m$
is the order of the differential equation.
Although this can be improved by using more complex operators the main target of 
the paper is the higher-order basis functions, hence this approach has not been pursued. 

Concerning the computational time, there is a clear optimum choice of $N_r$. On the one hand a small value of $N_r$ gives 
little computational cost but it is not very effective. On the other hand, large values of $N_r$ give more optimal multigrid
performance but at great expense per iteration.  
The optimal choice adopted in the remainder of this paper are $N_r=2,4,6,7,8,9$ for $N=4,8,12,16,20,24$, respectively.
This choice is not very sensitive as long as $N_r$ is not chosen on the low side where the computational cost rises quickly
with decreasing $N_r$.
\subsection{Final results on cubic domain}
\label{subsec:cubic_domain}
\begin{figure}[!ht]
       \centering
       \begin{subfigure}[b]{0.49\textwidth}
       \includegraphics[width=0.9\columnwidth]{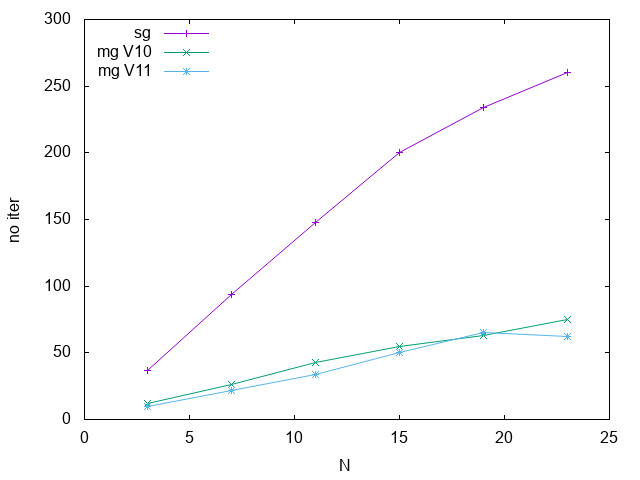}
       \caption{Iteration counts for constant angular basis functions.}
       \end{subfigure}
       \begin{subfigure}[b]{0.49\textwidth}
       \includegraphics[width=0.9\columnwidth]{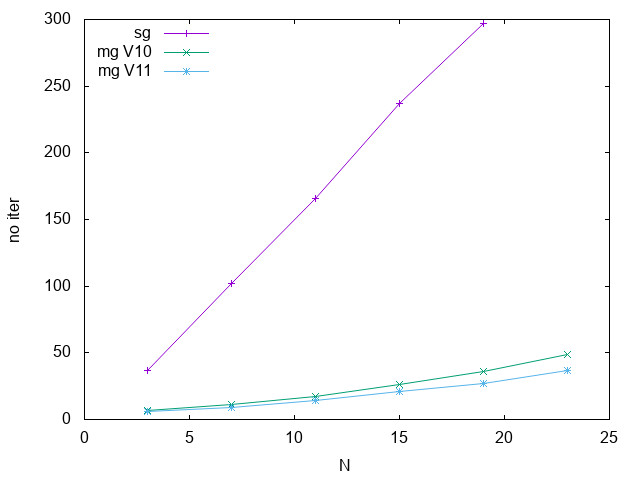}
       \caption{Iteration counts for linear angular basis functions.}
       \end{subfigure} \\
       \begin{subfigure}[b]{0.49\textwidth}
       \includegraphics[width=0.9\columnwidth]{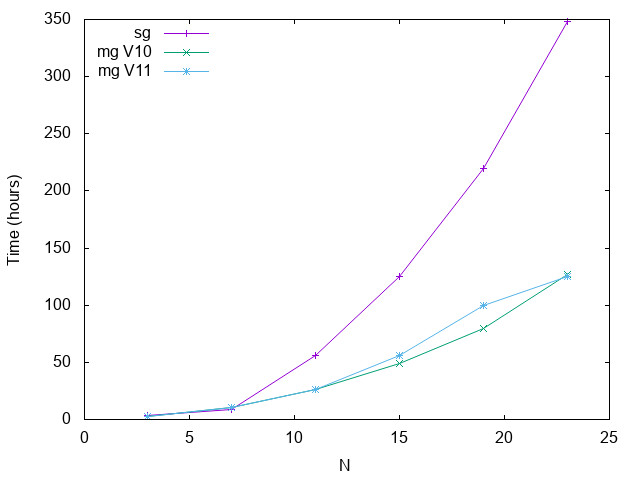}
       \caption{Computational time for constant angular basis functions.}
       \end{subfigure}
       \begin{subfigure}[b]{0.49\textwidth}
       \includegraphics[width=0.9\columnwidth]{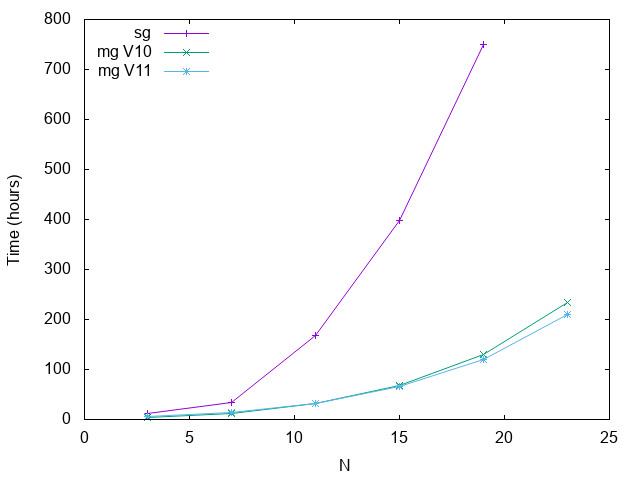}
       \caption{Computational time for linear angular basis functions.}
       \end{subfigure}
       \caption{Number of Krylov iterations and computational time for the standard sweep and
                the multigrid preconditioner using the V(1,0) and V(1,1) cycles for the 3D box geometry with uniform source and 
                isotropic angular refinement. Iteration counts given in a and b. Computational times are given in c and d. Constant 
                angular basis functions used in a and c. Linear basis is used in b and d.}
       \label{fig:cubic_final_results}
\end{figure}

To illustrate the efficiency of the final multigrid procedure using the alternative sweep on the coarsest grid as solver and 
the reduced scatter order in the multigrid preconditioner we compare the single mesh algorithm to the multigrid algorithm
in terms of the number of Krylov iterations and the computational time used.
Figure~\ref{fig:cubic_final_results} illustrates the main results where the parameters varied are the multigrid cycle used (V(1,0)
and V(1,1)) and the scatter order $N$ which is varied between the relatively smooth case of $N=4$ and very forward peaked $N=24$.
It is clear that the multigrid preconditioner gives a great enhancement of the performance compared to the single grid case
both in terms of the number of Krylov iterations required and more importantly in terms of computational effort. As seen before,
the constant basis functions do not lead to the same performance enhancements as the linear basis.
Both sets of basis functions show the multigrid advantage 
to increase as the scatter order increases. The computational time saved is between 3 for the constant basis case (at $N=24$) 
and 6.5 for the linear basis case (at $N=20$, as the single grid at $N=24$ were considered too expensive to run).
\subsection{Cubic domain with a boundary source}
\label{subsec:bc_cubic_domain}
\begin{figure}[!ht]
        \centering
        \begin{subfigure}[b]{0.32\textwidth}
        \includegraphics[width=0.9\columnwidth]{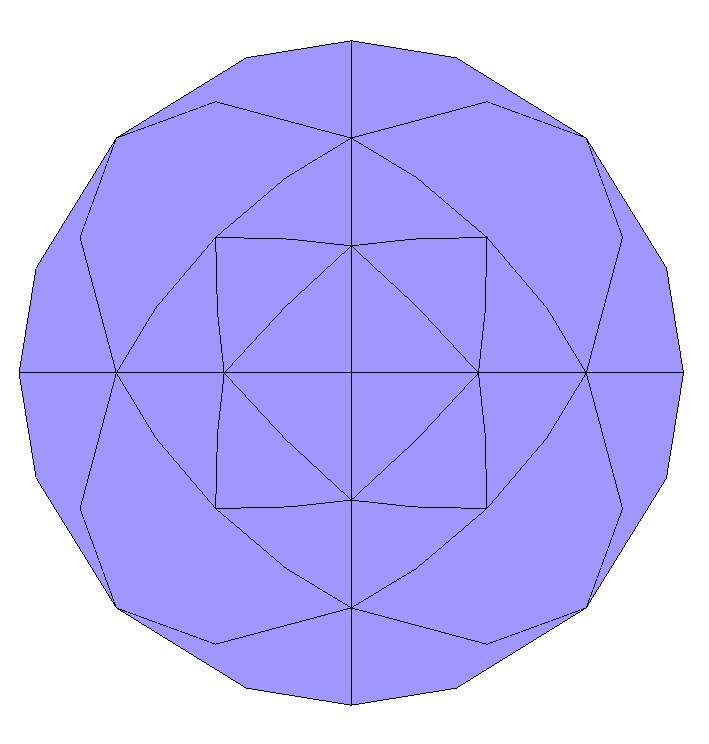}
        \end{subfigure}
        \begin{subfigure}[b]{0.32\textwidth}
        \includegraphics[width=0.9\columnwidth]{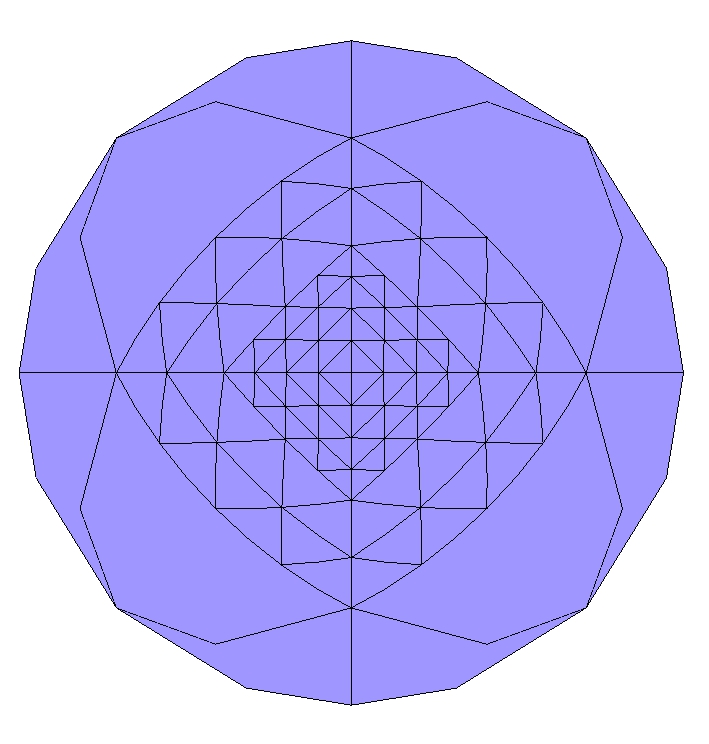}
        \end{subfigure}
        \begin{subfigure}[b]{0.32\textwidth}
        \includegraphics[width=0.9\columnwidth]{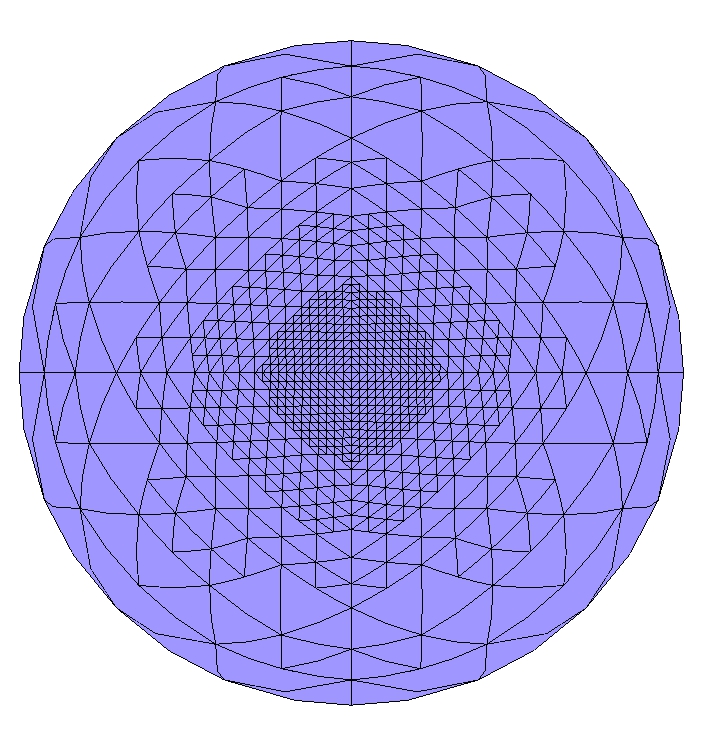}
        \end{subfigure}
        \caption{Non-isotropically refined angular meshes corresponding to $l_{max}$ = 2 (left), 4 (middle), 6 (right). 
                 The meshes are viewed from the positive $\Omega_z$-direction.}
        \label{fig:non_iso_meshes}
\end{figure}

\begin{figure}[!ht]
       \centering
       \begin{subfigure}[b]{0.49\textwidth}
       \includegraphics[width=0.9\columnwidth]{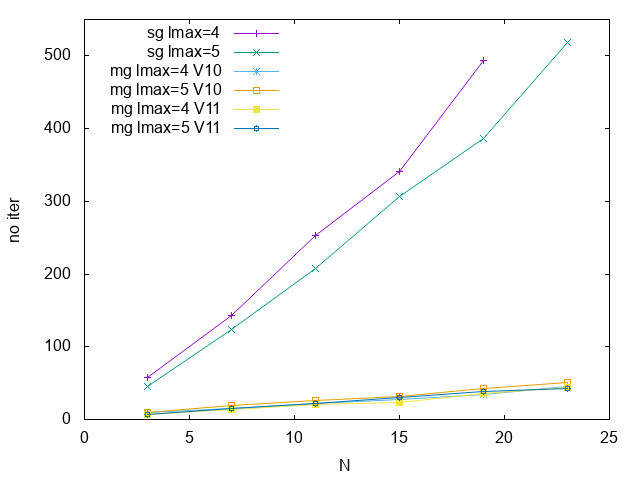}
       \caption{Iteration counts for constant angular basis functions.}
       \end{subfigure}
       \begin{subfigure}[b]{0.49\textwidth}
       \includegraphics[width=0.9\columnwidth]{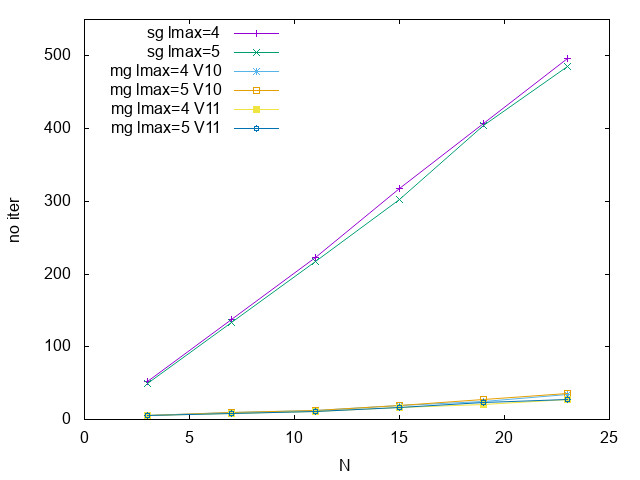}
       \caption{Iteration counts for linear angular basis functions.}
       \end{subfigure} \\
       \begin{subfigure}[b]{0.49\textwidth}
       \includegraphics[width=0.9\columnwidth]{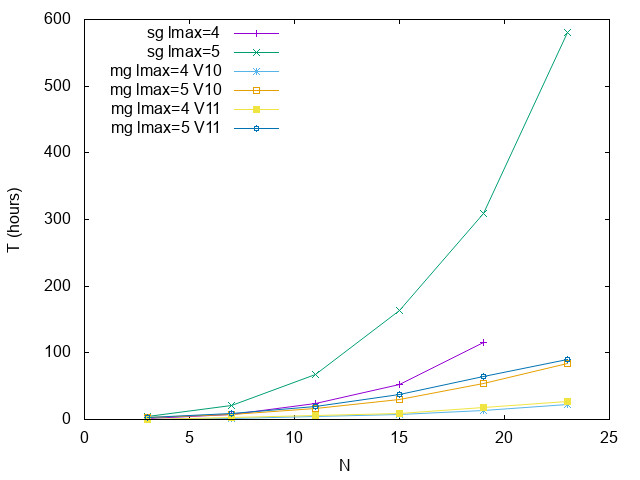}
       \caption{Computational time for constant angular basis functions.}
       \end{subfigure}
       \begin{subfigure}[b]{0.49\textwidth}
       \includegraphics[width=0.9\columnwidth]{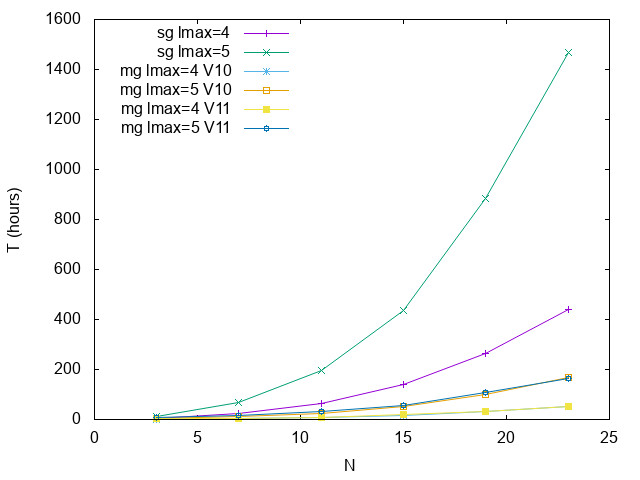}
       \caption{Computational time for linear angular basis functions.}
       \end{subfigure}
       \caption{Number of Krylov iterations and computational time for the standard sweep and
                the multigrid preconditioner using the V(1,0) and V(1,1) cycles for the 3D box geometry with boundary source and 
                anisotropic angular refinement with varying depth $l_{max}$. Iteration counts given in a and b. Computational times 
                are given in c and d. Constant 
                angular basis functions used in a and c. Linear basis is used in b and d.}
       \label{fig:bc_results}
\end{figure}

Following Turcksin and Morel \cite{Turcksin2012} and our previous work on Fokker-Planck solution using multigrid 
\cite{Lathouwers2019} we consider the same cubic domain as earlier but impose a boundary source composed of a 
pencil beam in the z-direction on part of the cube surface
\begin{equation}
  \phi (\boldr,\bOmega) = \delta(\boldr,\bOmega - \boldr,\bOmega_z), \; z=0, \; 2 < x,y < 3
\end{equation}
The Dirac function is used in the DG formalism (see \cite{Lathouwers2019}).
To capture the forward peaked radiative field due to the boundary condition used, we use a refined angular mesh
concentrating elements on the z-pole of the sphere. For a given maximum level $l_{max}$ in the angular, elements 
with $\Omega_z > 0.97$ are assigned this refinement level, elements with $0.9 < \Omega_z < 0.97$ are assigned level $l_{max} - 1$,
elements with $0.8 < \Omega_z < 0.9$ are assigned level $l_{max} - 2$, 
elements with $0.6 < \Omega_z < 0.8$ are assigned level $l_{max} - 3$, and all remaining elements are assigned level $l_{max} - 4$.
The angular meshes contain 20, 32, 44, 128, 440 and 1388 angular elements for $l_{max}$ between 1 and 6. 
Some of theses meshes are shown in Figure~\ref{fig:non_iso_meshes}. The non-isotropic refinement is one of the advantages 
of the present method over discrete ordinates techniques for highly non-isotropic radiative fields.
The results of the multigrid preconditioned method versus the regular sweep preconditioner are shown in 
Figure~\ref{fig:bc_results}.
As in the isotropically refined case, the multigrid based solver gets increasingly efficient w.r.t. to the single grid case
with higher scatter order $N$. Here, the savings are even greater, i.e. up to around a factor 
of 6 for the constant basis functions ($N=24$) and up to a factor of around 9 for the linear basis ($N=24$).
Finally, the effect of choice of multigrid cycling strategy is of little influence for the efficiency.
\section{Conclusions}
\label{sec:Conclusions}
The standard transport sweep is not an effective preconditioner for the transport equation in the presence of highly forward scatter. 
Diffusion synthetic acceleration is also known not to be effective in these cases. In this work we have adapted an angular multigrid preconditioner originally developed for the Fokker-Planck equation to the case using Legendre scatter modeling. Interpolation operators are defined naturally through the hierarchic nature of the discontinuous Galerkin angular discretization. Smoothing on each level is effected by the standard source iteration technique. A new sweep algorithm was developed that is better capable of solving the coarse mesh problem, even for high scatter orders.

The multigrid preconditioner was used in several test cases. The first test case consists of a 3D geometry with a uniform isotropic source. Uniform angular refinement is used. It was found that the most effective strategy is to use a reduced scatter order for the multigrid preconditioner. The optimal value
was found to increase slightly with the scatter order, $N$. Using this optimal value, a comparison was done with the standard sweep preconditioner in terms of number of Krylov iterations required and the computational time to solve the problem. In all cases the multigrid preconditioner outperformed the sweep algorithm, especially for the linear angular basis. With increasing scatter order, the difference 
between sweep and multigrid become greater. Another test case is a 3D geometry with a uni-directional boundary source. Here, 
the angular meshes are anisotropically refined to capture the anisotropic radiation field induced by the boundary condition. 
The multigrid preconditioner was again highly effective compared to the sweep preconditioner with similar trends as found in the volumetric source case.
\bibliography{./library}

\begin{thebibliography}{14}
\providecommand{\natexlab}[1]{#1}
\providecommand{\url}[1]{\texttt{#1}}
\expandafter\ifx\csname urlstyle\endcsname\relax
  \providecommand{\doi}[1]{doi: #1}\else
  \providecommand{\doi}{doi: \begingroup \urlstyle{rm}\Url}\fi

\bibitem[de~Oliveira et~al.(2000)de~Oliveira, Pain, and Eaton]{Oliveira2000}
C.R.E. de~Oliveira, C.C. Pain, and M.~Eaton.
\newblock {Hierarchical Angular Preconditioning for the Finite Element
  Spherical Harmonics Radiation Transport Method}.
\newblock In \emph{Proceedings of PHYSOR 2000, International Topical Meeting on
  Advances in Reactor Physics and Mathematics and Computation into the Next
  Millenium}. Americal Nuclear Society, 2000.

\bibitem[Drumm and Fan(2017)]{Drumm2017}
C.R. Drumm and W.C. Fan.
\newblock {Multilevel Acceleration of Scattering-source Iterations with
  Application to Electron Transport}.
\newblock \emph{Nuclear Engineering and Technology}, 49:\penalty0 1114--1124,
  2017.

\bibitem[Geuzaine and Remacle(2009)]{gmsh2009}
C.~Geuzaine and J.F. Remacle.
\newblock {Gmsh: A Three-dimensional Finite Element Mesh Gemerator with
  Built-in Pre- and Post-processing Facilities}.
\newblock \emph{Int. Journal for Numerical Methods in Engineering}, 79\penalty0
  (11):\penalty0 1309--1331, 2009.

\bibitem[Hennink and Lathouwers(2017)]{Hennink2017}
A.~Hennink and D.~Lathouwers.
\newblock {A Discontinuous Galerkin Method for the Mono-Energetic Fokker-Planck
  Equation based on a Spherical Interior Penalty Formulation}.
\newblock \emph{Journal of Computational and Applied Mathematics},
  330:\penalty0 253--267, 2017.

\bibitem[K{\'{o}}ph{\'{a}}zi and Lathouwers(2015)]{Kophazi2015}
J.~K{\'{o}}ph{\'{a}}zi and D.~Lathouwers.
\newblock {A Space-Angle DGFEM Approach for the Boltzmann Radiation Transport
  Equation with Local Angular Refinement}.
\newblock \emph{Journal of Computational Physics}, 297:\penalty0 637--668,
  2015.

\bibitem[Lathouwers and Perk\'o(2019)]{Lathouwers2019}
D.~Lathouwers and Z.~Perk\'o.
\newblock {An Angular Multigrid Preconditioner for the Radiation Transport
  Equation with Fokker-Planck Scattering}.
\newblock \emph{Journal of Computational and Applied Mathematics},
  350:\penalty0 165--177, 2019.

\bibitem[Lee(2010{\natexlab{a}})]{Lee2010a}
B.~Lee.
\newblock {A Novel Multigrid Method for SN Discretizations of the
  Mono-Energetic Boltzmann Transport Equation in the Optically Thick and Thin
  Regimes with Anisotropic Scatter, Part I}.
\newblock \emph{SIAM J. Sci. Comput.}, 31\penalty0 (6):\penalty0 4744--4773,
  2010{\natexlab{a}}.

\bibitem[Lee(2010{\natexlab{b}})]{Lee2010b}
B.~Lee.
\newblock {Improved Multiple-Coarsening Methods for SN Discretizations of the
  Bolltzmann Equation}.
\newblock \emph{SIAM J. Sci. Comput.}, 32\penalty0 (5):\penalty0 2497--2522,
  2010{\natexlab{b}}.

\bibitem[Morel and Manteuffel(1999)]{Morel1991}
J.E. Morel and T.A. Manteuffel.
\newblock {Two-Level Fourier Analysis of a Multigrid Approach for Discontinuous
  Galerkin Discretizations}.
\newblock \emph{Nucl. Sci. and Engng.}, 107:\penalty0 330--342, 1999.

\bibitem[Pain et~al.(2006)Pain, Eaton, Smedley-Stevenson, Goddard, Piggott, and
  de~Oliveira]{Pain2006}
C.C. Pain, M.D. Eaton, R.P. Smedley-Stevenson, A.J.H. Goddard, M.D. Piggott,
  and C.R.E. de~Oliveira.
\newblock {Space-Time Streamline Upwind Petrov-Galerkin Methods for the
  Boltzmann Transport Equation}.
\newblock \emph{Comput. Methods Appl. Mech. Engrg.}, 195:\penalty0 4334--4357,
  2006.

\bibitem[Pautz et~al.(1999)Pautz, Morel, and Adams]{Pautz1999}
S.~Pautz, J.E. Morel, and M.L. Adams.
\newblock {An Angular Multigrid Aceleration Method for the SN Equations with
  Highly Forward-Peaked Scattering}.
\newblock In \emph{Proceedings of International Conference on Mathematics and
  Computation, Reactor Physics and Environmental Analyses in Nuclear
  Applications}. Americal Nuclear Society, 1999.

\bibitem[Turcksin et~al.(2012)Turcksin, Ragusa, and Morel]{Turcksin2012}
B.~Turcksin, J.C. Ragusa, and J.E. Morel.
\newblock {Angular Multigrid Preconditioner for Krylov-Based Solution
  Techniques Applied to the {S}n Equations with Highly Forward-Peaked
  Scattering}.
\newblock \emph{Transport Theory and Statistical Physics}, 41\penalty0
  (1-2):\penalty0 1--22, 2012.

\bibitem[Van~der Vorst(1992)]{bicgstab}
H.A. Van~der Vorst.
\newblock {Bi-CGSTAB: A Fast and Smoothly Converging Variant of Bi-CG for the
  Solution of Nonsymmetric Linear Systems}.
\newblock \emph{SIAM J. Sci. and Stat. Comput.}, 13\penalty0 (2):\penalty0
  631--644, 1992.

\bibitem[Wesseling(1992)]{Wesseling1992}
P.~Wesseling.
\newblock \emph{{An Introduction to Multigrid Methods}}.
\newblock John Wiley and Sons, 1992.

\end{thebibliography}
\end{document}